\documentclass[a4paper,11pt,notitlepage]{amsart}
\usepackage{amsmath,amssymb,mathrsfs,amsthm}

\usepackage{verbatim}
\usepackage[spanish,USenglish]{babel} 
\usepackage{color}
\usepackage{graphicx}

\newtheorem{thm}{Theorem}[section]
\newtheorem{corollary}[thm]{Corollary}

\newtheorem{theorem}[thm]{Theorem}
\newtheorem{lemma}[thm]{Lemma}

\theoremstyle{definition}

\theoremstyle{remark}
\newtheorem{remark}[thm]{Remark}

\numberwithin{equation}{section}

\newcommand\Real{{\mathfrak R}{\mathfrak e}\,} %

\newcommand{\R}{\mathbb{{R}}}
\newcommand{\Z}{\mathbb{{Z}}}
\newcommand{\N}{\mathbb{{N}}}
\newcommand{\T}{\mathbb{{T}}}
\newcommand{\C}{\mathbb{{C}}}

\newcommand{\md}{\mathrm{{d}}}
\begin{document}

\title[Large time behaviour for the heat equation on $\Z$]{Large time behaviour for the heat equation on $\Z,$ moments and decay rates}

\author[L. Abadias]{Luciano Abadias}
\address[L. Abadias]{Departamento de Matem\'aticas, Instituto Universitario de Matem\'aticas y Aplicaciones, Universidad de Zaragoza, 50009 Zaragoza, Spain.}
\email{labadias@unizar.es}

\author[J. González-Camus]{Jorge González-Camus}
\address[J. González-Camus]{Departamento de Matemáticas y Ciencias de la Computación, Facultad de
Ciencia, Universidad de Santiago de Chile, Santiago, Chile.}
\email{jorge.gonzalezcam@usach.cl}

\author[P. J. Miana]{Pedro J. Miana}
\address[P. J. Miana]{Departamento de Matem\'aticas, Instituto Universitario de Matem\'aticas y Aplicaciones, Universidad de Zaragoza, 50009 Zaragoza, Spain.}
\email{pjmiana@unizar.es}

\author[J. C. Pozo]{Juan C. Pozo}
\address[J. C. Pozo]{Departamento de Matemáticas, Facultad de Ciencias, Universidad de Chile, Las
Palmeras 3425, \~Nu\~noa, Santiago, Chile.}
\email{jpozo@ug.uchile.cl}

\thanks{The first and third named authors have been partly supported by Project MTM2016-77710-P, DGI-FEDER, of the MCYTS and Project E26-17R, D.G. Arag\'on, Universidad de Zaragoza, Spain. The first author has been also supported by Project for Young Researchers, Fundación Ibercaja and Universidad de Zaragoza, Spain. The second author has been supported by the National Agency for Research and Development (ANID),
Beca Doctorado Nacional Chile 2016, Number 21160120. The fourth author has been partially supported by Chilean research Grant FONDECYT grant 1181084.}



\begin{abstract}
The paper is devoted to understand the large time behaviour and decay of the solution of the discrete heat equation in the one dimensional mesh $\Z$ on $\ell^p$ spaces, and its analogies with the continuous-space case. We do a deep study of the moments of the discrete gaussian kernel (which is given in terms of Bessel functions), in particular the mass conservation principle; that is reflected on the large time behaviour of solutions. We prove asymptotic pointwise and $\ell^p$ decay results for the fundamental solution. We use that estimates to get rates on the $\ell^p$ decay and large time behaviour of solutions. For the  $\ell^2$ case, we get optimal decay by use of Fourier techniques.
\end{abstract}

\date{}

\maketitle

\section{Introduction}

Diffusion processes happen constantly around us. Such diffusions processes have been studied along the history by many researchers of several scientific branches. For that purpose, many mathematical models have been appeared in the literature in the last centuries. One of the main mathematical areas which studies these processes is the mathematical analysis, using techniques from PDEs, functional analysis, harmonic analysis, among others. In particular, parabolic equations can be used to know the behaviour of diffusion for certain ``objects" or ``things'' of distinct nature, as heat, diseases, population growth, and fashions. The dynamic of the evolution models is one of the main feature studied by the mathematician because it contributes to know, in particular, the large time behaviour of the diffusion.\\

The nature of the object that one wants to study defines the framework of the model; for example, it seems natural that the model will be different if we want to understand the diffusion of a type of particles in a liquid, than if we want to know the evolution of the diffusion of an idea in the society. In the first case, the space where the diffusion takes places could be $\R^3,$ or a subset of it; in the second one case would be a network.\\

The study of parabolic models usually first looks at to the most known classical diffusion equation, the heat equation, which was introduced by J. Fourier in 1822, see \cite{Fourier}. The heat equation is governed by the Laplace operator $\Delta;$ such differential operator takes several forms depending on the spaces where we are working.  In this manuscript, we focus on the heat propagation model on the infinite one dimensional mesh $\Z$ (also called one dimensional infinite lattice). In the following, we will refer to the operator who defines such model as the discrete Laplacian, which is given by $$\Delta_\md f(n):=f(n+1)-2f(n)+f(n-1),\quad n\in\Z,$$ for each sequence $f$ defined on $\Z.$ The main goal of the paper is to know the decay rates and large time behaviour on $\ell^p:=\ell^p(\Z)$ ($1\leq p\leq \infty$) for the solution of the discrete heat problem  \begin{equation}\label{eq0}
\left\{\begin{array}{ll}
\partial_t u(t,n)=\Delta_\md u(t,n)+g(t,n),&n\in\Z,\,t>0,\\ \\
u(0,n)=f(n),&n\in\Z.
\end{array} \right.
\end{equation}
In the previous equation, $u$ denotes the solution, $f$ the initial data, and $g$ the linear forcing term. It is known that if $f\in\ell^p,$ and $g(t)\equiv g(t,\cdot)$ belongs to $L^1_{loc}([0,\infty),\ell^p),$ then \eqref{eq0} has a unique mild solution (see \cite[Part I, Chapter 3]{ABHN} and \cite[Chapter VI, Section 7]{EN}) given by \begin{equation}\label{sol}u(t,n)=W_t f(n)+\int_0^{t}(W_{t-s}g(s,\cdot))(n)\,ds,\end{equation} where $(W_t)_{t\geq 0}$ is the heat semigroup whose infinitesimal generator is $\Delta_{\md};$ more explicitly, the semigroup $(W_t)_{t\geq 0}$ is of convolution type, that is $$W_t f(n)=(G(t,\cdot)*f)(n):=\sum_{j\in\Z}G(t,n-j)f(j),\quad n\in\Z,$$ where the discrete heat kernel is given by  $G(t,n):=e^{-2t}I_{n}(2t),$ being $I_n$ the Bessel function of imaginary argument and order $n\in\Z.$ Moreover, if $g\in C([0,\infty),\ell^p)$ then $u$ given by $\eqref{sol}$ is a classical solution of \eqref{eq0} (see \cite[Part I, Chapter 3]{ABHN} and \cite[Chapter VI, Section 7]{EN}).\\

The convolution discrete heat kernel is called the fundamental solution of our diffusion process; it is the solution of \begin{equation}\label{eq1}
\left\{\begin{array}{ll}
\partial_t u(t,n)=\Delta_\md u(t,n),&n\in\Z,\,t> 0,\\ \\
u(0,n)=\delta_{0}(n),&n\in\Z,
\end{array} \right.
\end{equation}
where the initial data is the unit impulse symbol $\delta_0$ (the Dirac mass on $\Z$). Recall that $\delta_0(0)=1$ and $\delta_0(n)=0$ for $n\neq 0.$\\

Under our knowledge, H. Bateman was the first author who proposed the solution of \eqref{eq1} in \cite{Bateman}. Moreover, he studied a broad set of differential-difference equations (heat and wave equations), whose solutions are given in terms of special functions; the Bessel function $J_n,$ the Bessel function of imaginary argument $I_n,$ the Hermite polynomial $H_n$ and the exponential function. After that, several papers on the topic have appeared. For example, the authors in \cite{C} do an harmonic  study of the previous problem, studying maximum principles for the discrete fractional Laplacian, weighted $\ell^p$-boundedness of conjugate harmonic functions, Riesz transforms and square functions of Littlewood-Paley. They also prove that $(W_t)_{t\geq 0}$ is a Markovian $C_0$-semigroup on $\ell^p$ whose generator is $\Delta_\md$ (observe that $\Delta_\md$ is a bounded operator on $\ell^p$). This means that it has the semigroup property ($W_tW_s=W_{t+s}$), the strong continuous property ($f=W_0f=\lim_{t\to0^+}W_t f$ in $\ell^p$), it is contractive ($\|W_t f\|_p\leq \|f\|_p$), it is positive ($W_t f\geq 0$ whenever $f\geq 0$), and it satisfies the mass conservation principle ($W_t1=1$). In \cite{C2}, the authors obtained a series convolution representation for the fractional powers of the discrete Laplacian, and they used it to prove regularity properties on discrete H\"older spaces and convergence results from the discrete to the continuous case. In \cite{LR}, a deep study for the discrete wave problem is done. The authors also study spectral properties on $\ell^p$ spaces. They get that $(W_t)_{t\geq 0}$ is an holomorphic semigroup on $\ell^p$ of angle $\pi,$ with the spectrum $\sigma(\Delta_{\md})=[-4,0].$ Recently, some super-diffusive linear processes in the one-dimensional mesh are proposed in \cite{Estrada} by means of $d$-Laplacians.\\

As we have said, the main aim in this paper is to study the large time behaviour of solutions of \eqref{eq0}, and also to compare the results obtained to the continuous case. Recall that the solution of the heat equation on $L^p(\R)$ \begin{equation}\label{eq6}
\left\{\begin{array}{ll}
\partial_t v(t,x)=\Delta v(t,x),&x\in\R,\,t> 0,\\ \\
v(0,x)=f(x),&x\in\R,
\end{array} \right.
\end{equation}
is $v(t,x)=\int_{\R}g_t(x-y)f(y)\,dy$ (assuming $f\in L^p(\R)$), where $g_{t}(x):=\frac{1}{\sqrt{4\pi t}}e^{-\frac{x^2}{4t}}$ is the gaussian kernel and $\Delta=\frac{\partial^2}{\partial x^2}$. Observe that integrating over all of $\R,$ we get that the total mass of solutions is conserved for all time, that is, $$\int_{\R}v(t,x)\,dx=\int_{\R}v(0,x)\,dx=\int_{\R}f(x)\,dx.$$ This fact is known as the mass conservation principle. It is known that the time invariance property is reflected in the large time behaviour of solutions; if $M=\int_{\R}f(x)\,dx,$ then \begin{equation}\label{0.1}t^{\frac{1}{2}(1-1/p)}\|v(t)-M g_{t}\|_{p}\to0,\quad\text{as }\,t\to \infty,\end{equation} for $1\leq p\leq \infty,$ that is, the rate of $L^p$-convergence from $v$ to $M g_t$ is $o\biggl(\frac{1}{t^{\frac{1}{2}(1-1/p)}}\biggr)$ as $t$ goes to infinity. \\

On the other hand, the $p$-energy for $p>1$ is not conservative. It is known that  $$\|v(t,\cdot)\|_{p}\leq C \|f\|_{q}t^{-\frac{1}{2}(1/q-1/p)},\quad \|\nabla v(t,\cdot)\|_{p}\leq C \|f\|_{q}t^{-\frac{1}{2}(1/q-1/p)-\frac{1}{2}},$$ $$\|\partial_t v(t,\cdot)\|_{p}\leq C \|f\|_{q}t^{-\frac{1}{2}(1/q-1/p)-1},$$ for $f\in L^q(\R)$ and $1\leq q\leq p\leq\infty.$ Such inequalities can be generalized to higher derivative orders, see \cite[Part II, Section 9.3]{ebert}.\\

Also, one can study the first moment, the vector quantity $\int_{\R}x\,v(t,x)\,dx.$ If we assume that $(1+|x|)f\in L^1(\R),$ such moment is also conserved in time. This also allows to get an improvement on the convergence rate \eqref{0.1}; if $(1+|x|)f\in L^1(\R),$ then \begin{equation}\label{0.2}t^{\frac{1}{2}(1-1/p)}\|v(t)-M g_{t}\|_{L^p(\R)}\leq Ct^{-1/2}.\end{equation} However, the second moment $$\int_{\R}x^2 \,v(t,x)\,dx=\int_{\R}x^2\,f(x)\,dx+2t\int_{\R}\,f(x)\,dx$$ is time dependent. Moreover, only integral quantities conserved by the solutions of \eqref{eq6} are the mass and the first moment. Some references where one can find previous known asymptotic results about the solution of \eqref{eq6} are \cite{DZ,ZE,F,V,Z}.\\

This type of large-time asymptotic results have been also studied for several diffusion problems. For example in \cite{DD,ZE,GV,KS,N} the authors studied large-time behaviour and other asymptotic estimates for the solutions of different diffusion problems in $\R^N,$ and similar aspects are studied for open bounded domains in \cite{D2, GV2}. Estimates for heat kernels on manifolds have been also studied in \cite{D,Gr,L}. In \cite{M}, the author obtains gaussian upper estimates for the heat kernel associated to the sub-Laplacian on a Lie group, and also for its first-order time and space derivatives.\\

Along the paper, we revisit the results commented previously, but now in the discrete setting. First, in Section 2, we state some known inequalities and identities, mainly for the Bessel function of imaginary argument $I_n,$ which will be applied in next sections. The Section 3 is devoted to the moments of Bessel functions $I_n,$ and in particular of the fundamental solution $G(t,n).$ It allows to get the mass conservation principle, and the first and second order moments to the solution of \eqref{eq0} with null linear forcing term ($g\equiv 0$). This suggests that it is possible to obtain similar asymptotic results to the continuous case. In Section 4 we do a deep study for the large time behaviour of the discrete gaussian kernel $G(t,n).$ We get a technical result (Lemma \ref{Pointwiseestimates}) where we study pointwise estimates, which will be a key role on the main result of the paper, and we also study $\ell^p$ decay rates and compare to the continuous case ones. The following section, Section 5, is focused on studying $\ell^p$ decayment for the solutions of $\eqref{eq0}$ assuming $\ell^q$ conditions ($1\leq q\leq p\leq \infty$) on the data $f$ and on the linear forcing term $g.$ In addition, we get an optimal $\ell^2$ result by Fourier techniques when the data belongs to $\ell^1.$ In Section 6 we present how the solution of \eqref{eq0} converges, for large $t,$ to a constant (depending on the mass of the problem) times the discrete gaussian, on $\ell^p$ spaces, similarly to what happens in the continuous case (\eqref{0.1} and \eqref{0.2}). Finally, in the last section (Section 7), we do some remarks and pose open questions.

\section{Preliminaries and known results}

In this section we show some known identities, inequalities and results that we will use along the paper, to make easier the reading.

First, we introduce the following inequality \begin{equation}\label{eq2.2}
(1-r)^{\eta}r^{\gamma}\leq \biggl(\frac{\gamma}{\gamma+\eta}\biggr)^{\gamma},\quad \gamma\geq 0,\eta>0,0<r<1,
\end{equation}
which was a key point in the proof of many results in \cite{C}.

Let $\alpha,z\in\C,$ then the gamma function satisfies $$\frac{\Gamma(z+\alpha)}{\Gamma(z)}=z^{\alpha}\biggl(1+\frac{\alpha(\alpha+1)}{2z}+O(|z|^{-2})\biggr),\quad |z|\to\infty,$$ whenever $z\neq 0,-1,-2,\ldots$ and $z\neq -\alpha,-\alpha-1,\ldots,$ see \cite[Eq.(1)]{ET}. Particularly

\begin{equation}\label{double3}
\frac{\Gamma(z+\alpha)}{\Gamma(z)}=z^{\alpha}\biggl(1+O\biggl({1\over |z|}\biggr)\biggr), \quad z\in\C_+,\,\Real\alpha>0.
\end{equation}

We denote by $I_n$ the Bessel function of imaginary argument (also called modified Bessel function of first kind) and order $n\in\Z,$ given by \begin{equation*}\label{DefBessel}I_n(t)=\sum_{m=0}^{\infty}\frac{1}{m!\Gamma(m+n+1)}\biggl(\frac{t}{2}\biggr)^{2m+n},\quad n\in\Z, z\in \C.\end{equation*}
Since $n$ is an integer, ${1\over \Gamma(n)}$ is taken to equal zero if $n=0, -1,-2,\ldots,$ so $I_n$ is defined on the whole complex plane, being an entire function. Now we give some known properties about Bessel functions $I_n$  which can be found in \cite[Chapter 5]{Leb} and \cite{W}. They satisfy that $I_{-n}=I_n$ for $n\in\Z$, $I_0(0)=1$, $I_n(0)=0$ for $n\not =0,$ and $I_n(t)\geq 0$ for  $n\in\Z$ and $t\geq0$. Also, the function $I_n$ has the semigroup property (also called Neumann's identity) for the convolution on $\Z,$ that is, \begin{equation*}\label{Semigroup}I_n(t+s)=\sum_{m\in\Z}I_m(t)I_{n-m}(s)=\sum_{m\in\Z}I_m(t)I_{m-n}(s),\quad t,s\geq 0,\end{equation*} see \cite[Chapter II]{Feller}. The generating function of the Bessel function $I_n$ is given by \begin{equation}\label{generating}e^{\frac{t(x+x^{-1})}{2}}=\sum_{n\in\Z}x^n I_n(t), \qquad x\neq 0, \,t\in \C.\end{equation} It also satisfies the following differential-difference equation
\begin{equation}\label{prop3}\frac{\partial}{\partial_t}I_n(t)=\frac{1}{2}\biggl(I_{n-1}(t)+I_{n+1}(t)\biggr), \quad z\in \C.\end{equation}

In the following we enumerate some integral representations of the modified Bessel functions which are useful along the paper;  let $n\in\N_0:=\N\cup\{0\}$ and  $t\in\C,$ it follows
\begin{eqnarray}\label{prop4}I_n(t)&=&\frac{t^n}{\sqrt{\pi}2^n\Gamma(n+1/2)}\int_{-1}^1 e^{-ts}(1-s^2)^{n-1/2}\,ds,\\
\label{prop5}I_{n+1}(t)-I_{n}(t)&=&-\frac{t^n}{\sqrt{\pi}2^n\Gamma(n+1/2)}\int_{-1}^1 e^{-ts}(1+s)(1-s^2)^{n-1/2}\,ds,\\
\label{prop6}
\displaystyle I_{n+2}(t)-2I_{n+1}(t)+I_n(t)&=&\frac{t^n}{\sqrt{\pi}2^n\Gamma(n+1/2)}\times\nonumber\\
&&\displaystyle \biggl(\frac{2}{t}\int_{-1}^1 e^{-ts}s(1-s^2)^{n-1/2}\,ds+\int_{-1}^1 e^{-ts}(1+s)^2(1-s^2)^{n-1/2}\,ds\biggr).
\end{eqnarray}

The following estimates of the Bessel functions are also well-known,
\begin{equation}\label{prop2}I_n(t)= \frac{e^{t}}{(2\pi t)^{1/2}}\biggl( 1-\frac{(n^2-1/4)}{4t}+O\biggl(\frac{1}{t^2}\biggr) \biggr), \quad t \text{ large},\end{equation}

\begin{equation}\label{prop8}I_n(t)\sim \frac{t^n}{2^n n!}, \quad t\to 0,n\in\N_0.\end{equation}

\begin{remark}\label{dirac}{\it Note that the property $I_n(0)=0$ for $n\neq 0$ and $I_0(0)=1,$ can be also obtained from \eqref{prop8}. This implies that $G(t,\cdot)$ converges pointwise to the dirac mass $\delta_0$ (which is the identity on discrete convolution) as $t\to 0.$}
\end{remark}

From the theory of Confluent Hypergeometric Functions, see \cite[Section 9.11]{Leb}, we have \begin{equation}\label{Hyper}\int_0^1 e^{-4ts} s^{\gamma-\alpha-1}(1-s)^{\alpha-1}\,ds=\Gamma(\gamma-\alpha)e^{-4t}\sum_{k=0}^{\infty}\frac{(4t)^k}{k!}\frac{\Gamma(\alpha+k)}{\Gamma(\gamma+k)}, \end{equation}
and
\begin{equation}\label{lebedev}
I_n(t)={(2n)!\over n!\Gamma(n+{1\over 2})}\biggl({t\over 2}\biggr)^ne^{-t}\sum_{k=0}^\infty{\Gamma(n+k+{1\over 2})\over (2n+k)!}{(2t)^k\over k!}
\end{equation}
for $n\in\N_0$ and $t\in \C$, see  \cite[(9.13.14)]{Leb}.

Let $\T$ be the one dimensional torus. In the following we identify $\T$ with the interval $(-\pi ,\pi],$ and the functions defined on $\T$ with $2\pi$-periodic functions on $\R$. So, we consider the Fourier transform $$\mathcal{F}(f)(\theta)=\sum_{n\in\Z}e^{i n\theta}f(n),\quad \theta\in\T,f\in\ell^1. $$ The Fourier transform is extended to $\ell^2$ and it is an isometry from $\ell^2$ to $L^2(\T),$ with
inverse operator $$\mathcal{F}^{-1}(\varphi)(n)=\frac{1}{2\pi}\int_{\T}e^{-i n \theta} \varphi(\theta)d\,\theta.$$ In the proof of  \cite[Proposition 1]{C} it is shown that $$G(t,n)=\frac{e^{-2t}}{2\pi}\int_{\T}e^{-i n \theta}e^{2t\cos \theta}\,d\theta,$$ and therefore $$\mathcal{F}(G(t,\cdot))(\theta)=e^{-2t(1-\cos\theta)}=e^{-4t\sin^2\frac{\theta}{2}}.$$


\section{Moments of Bessel functions}

In the following we study the moments of the gaussian discrete kernel, and we apply that to get the moments of the solution of \eqref{eq0} with null linear forcing term ($g\equiv 0$). We compare the results to the continuous case.

\begin{lemma}\label{sequence} Let   $(p_k)_{k\ge 0}$ be the sequence of polynomials $(p_k)_{k\ge 0}$ given by $p_0(t):=1$ and
$$
\left\{
\begin{array}{ll}
p_k'(t):=\displaystyle{\sum_{j=0}^{k-1}{2k\choose 2j}p_j(t)},&\,\\
p_k(0):=0&\,
\end{array}
\right.
$$
for $k\ge 1$ and $t\in \C$. Then the polynomial $p_k$ has positive integer coefficients, its degree is $k$, and $p_k(t)> 0$ for $t> 0$.
For $k\ge 1$, we write
$$p_{k}(t)=\sum_{n=1}^k a_{k,n}t^n, \qquad t\in \C;
$$
then $a_{k,1}=1$, and
$$
a_{k,n}={1\over n}\sum_{j=n-1}^{k-1}{2k\choose 2j}a_{j,n-1}
$$
for $ 2\le n\le k$; in particular $a_{k,2}= 4^{k-1}-1 $, $a_{k, k}=(2k-1)!!$.
\end{lemma}

\begin{proof} By definition, it is clear that the coefficients of the polynomial $p_k$ are natural numbers, its degree is $k$, and $p_k(t)> 0$ for $t> 0$. For $k\ge 1$, we write $p_{k}(t)=\sum_{n=1}^k a_{k,n}t^n$ with $a_{k,n}\in \N$. Then we have
\begin{eqnarray*}
p_k'(t)&=&\sum_{j=0}^{k-1}{2k\choose 2j}p_j(t)=1+\sum_{j=1}^{k-1}{2k\choose 2j}\sum_{n=1}^j a_{j,n}t^n= 1+\sum_{n=1}^{k-1}t^n \sum_{j=n}^{k-1} {2k\choose 2j}a_{j,n}.\cr
\end{eqnarray*}
We integrate to obtain that
$$
p_k(t)= t+ \sum_{n=1}^{k-1}t^{n+1} {1\over n+1}\sum_{j=n}^{k-1} {2k\choose 2j}a_{j,n}= t+ \sum_{n=2}^{k}t^{n} {1\over n}\sum_{j=n-1}^{k-1} {2k\choose 2j}a_{j,n-1}.
$$
Now  for $n=2$, we have that
$$
a_{k,2}={1\over 2}\sum_{j=1}^{k-1}{2k\choose 2j}a_{j,1}={1\over 2}\sum_{j=0}^{k}{2k\choose 2j}-1=2^{2k-2}-1=4^k-1,
$$
for $k\ge 2$, where we have applied \cite[Formula 0.151(2)]{G}. Finally for $n=k$, by induction method one gets
$$
a_{k,k}={1\over k} {2k\choose 2(k-1)}a_{k-1,k-1}= {1\over k} {2k\choose 2(k-1)}(2k-3)!!=(2k-1)!!.
$$
\end{proof}

\begin{remark}{\it Here we present the polynomials $p_k$ for values $0\le k\le 6$.
\begin{eqnarray*}
p_0(t)&=&1,\cr
p_1(t)&=&t,\cr
p_2(t)&=&t+3t^2,\cr
p_3(t)&=&t+15t^2+15t^3,\cr
p_4(t)&=&t+63t^2+210t^3+105t^4,\cr
p_5(t)&=&t+225t^2+2205t^3+3150t^4+945t^5,\cr
p_6(t)&=&t+1023t^2+21120t^3+65835t^4+51975t^5+10395t^6.\cr
\end{eqnarray*}
}
\end{remark}

In the next theorem, we present the moments (of arbitrary order) of the modified Bessel functions $(I_n)_{n\in \Z}$.

\begin{theorem} \label{pepe} For $k\ge 0$, we have that
$$
\sum_{n\in\Z} n^{2k}I_n(t)= e^{t}p_k(t), \qquad \sum_{n\in\Z} n^{2k+1}I_n(t)= 0, \qquad t\in \C,
$$
where the polynomials $(p_k)_{k\ge 0}$ are given in Lemma \ref{sequence}.
\end{theorem}
\begin{proof} Since $I_{-n}=I_{n}$, it is straightforward to check that
$\displaystyle{\sum_{n\in\Z} n^{2k+1}I_n(t)= 0.}$  By the generating formula (\ref{generating}), we have that
$$
{\partial^j\over \partial x^j}\left(e^{{t(x+x^{-1})\over 2}}\right)=\sum_{n\in\Z}n(n-1)\cdots(n-j+1)x^{n-j} I_n(t), \qquad j\ge 1.
$$
We take the value $x=1$ on previous identity and conclude that
$$
 e^{t}Q(t)=\sum_{n\in\Z}n(n-1)\cdots(n-j+1) I_n(t),\qquad t\in \C,
$$
for a polynomial $Q$. Now we write
$$
e^{t}P_k(t)=\sum_{n\in\Z} n^{2k}I_n(t), \qquad t\in \C,
$$
for some polynomials $(P_k)_{k\ge 0}$. Note that $P_0(0)=1$ and $P_k(0)=0$ for $k\not=0$. We derive to get
$$
e^{t}(P_k(t)+P'_k(t))={1\over 2}\sum_{n\in\Z} n^{2k}(I_{n-1}(t)+I_{n+1})(t)= {1\over 2}\sum_{n\in\Z} ((n+1)^{2k}+(n-1)^{2k})I_{n}(t),
$$
for $t\in \C,$ where we have used \eqref{prop3}. Since  $\displaystyle{(n+1)^{2k}+(n-1)^{2k}=2\sum_{j=0}^k{2k \choose 2j}n^{2j}}$, we obtain that
$$
e^{t}P'_k(t)= \sum_{n\in\Z} \sum_{j=0}^{k-1}{2k\choose 2j}n^{2j}I_{n}(t)=e^{t}\sum_{j=0}^{k-1}{2k\choose 2j}P_j(t),
$$
for $t\in \C$. Therefore $P_k=p_k.$
\end{proof}

For the discrete heat kernel $(G(t,\cdot))_{t\in \C}$, we obtain the following moments.
\begin{corollary} \label{mean}For $k\ge 0$, we have that
$$
\sum_{n=-\infty}^\infty n^{2k}G(t,n)= p_k(2t), \qquad \sum_{n=-\infty}^\infty n^{2k+1}G(t,n)= 0, \qquad t\in \C,
$$
where the polynomials $(p_k)_{k\ge 0}$ are given in Lemma \ref{sequence}.

\end{corollary}

\begin{remark}{\it Note that for the continuous heat (or Gaussian) semigroup $(g_{t})_{\Real t>0}$ we have
$$
\int_{\R} s^{2k} g_{t}(s)\,ds={(2k)!\over k!}t^k, \qquad\int_{\R} s^{2k+1} g_{t}(s)\,ds=0, \qquad  \Real t>0,
$$
see \cite[Formula 3326]{G}. It is interesting to check that
$$
\displaystyle{\lim_{t\to + \infty}{\sum_{n\in\Z} n^{2k}G(t,n)\over \int_{\R} s^{2k} g_{t}(s)\,ds}=\lim_{t\to + \infty}{p_k(2t)\over {(2k)!\over k!}t^k}= {(2k-1)!!2^k\over {(2k)!\over k!}}}=1.
$$
}
\end{remark}

Finally, we study the moments of the solution of \eqref{eq0} with null linear forcing term ($g\equiv 0$).

\begin{remark}{\it From Corollary \ref{mean} it follows that $\sum_{n\in\Z}G(t,n)=1.$ This fact leads to observe the mass conservation principle for the diffusion problem \eqref{eq0} with $g\equiv 0$; $$\sum_{n\in\Z} f(n)=\sum_{n\in\Z} u(t,n),\quad t\geq 0.$$ We will see in Section \ref{Non-homogeneous Problem} that the mass conservation principle is reflected in the large time behaviour of the solution of \eqref{eq0}. Moreover, the first moment is also conservative; if $nf\in \ell^1$ one gets $$
\partial_t \sum_{n\in\Z}n\,u(t,n)=\sum_{n\in\Z}n \Delta_\md u(t,n)=\sum_{n\in\Z}((n-1)-2n+(n+1))u(t,n)=0,$$ therefore $$\sum_{n\in\Z}n \, u(t,n)=\sum_{n\in\Z}n\,f(n).$$ However, similarly to the heat problem on $\R,$ by Corollary \ref{mean} it follows that the second moment is time-dependent (whenever $n^2 f \in\ell^1$): \begin{eqnarray*}
\sum_{n\in\Z}n^2u(t,n)&=&\sum_{n\in\Z} n^2\sum_{j\in\Z}G(t,n-j)f(j)=\sum_{j\in\Z}f(j)\sum_{n\in\Z} ((n-j)^2+j^2+2j(n-j))G(t,n-j)\\
&=&\sum_{j\in\Z}j^2f(j)+2t\sum_{j\in\Z}f(j).
\end{eqnarray*}}
\end{remark}

\section{Pointwise and $\ell^p$ estimates of the fundamental solution}\label{point}

In this section we show pointwise and $\ell^p$ decay rates for $G(t,n).$ In the following we denote the first forward difference $\nabla_\md f(n):=f(n+1)-f(n),$ for $n\in\Z$ and $f:\Z\to \C.$ The results obtained along the paper where this first difference appears would be also right for the backward difference, $\nabla_{\md,-} f(n):=f(n)-f(n-1),$ for $n\in \Z.$ The notation $\nabla_\md$ and not $\nabla_{\md,+}$ is for convenience, because we will only use during whole the paper the forward operator $\nabla_{\md}.$ It is known that the operators $\nabla_{\md,\pm}$ generate markovian $C_0$-semigroups on $\ell^p,$ see \cite{ADT}.

The next result shows pointwise asymptotic estimates for $t$ large for the semigroup $G(t,n),$ depending on whether  $R:=\frac{|n|^2}{t}$ is greater or less than 1. Some items on next result do not depends on $R,$ but we do such a distinction because it is a key technical lemma that we will use in that form in the large time behaviour of solutions (Section \ref{Non-homogeneous Problem}).

\begin{lemma}\label{Pointwiseestimates}
Let $R=\frac{|n|^2}{t}$ with $n\in\Z\setminus\{0\}$ and $t>0$ large enough. Then there is $C>0$ such that
\begin{itemize}

\item[(i)] $|G(t,n)|\leq \frac{C}{t^{1/2}}$ for $R\leq 1,$ and $|G(t,n)|\leq \frac{C}{|n|^{3}}$ for $R\geq 1.$\\

\item[(ii)] If $n\in\N,$  $|\nabla_{\md} G(t,n)|\leq \frac{C|n|}{t^{3/2}}$ for $R\leq 1,$ and $|\nabla_{\md} G(t,n)|\leq \frac{Ct}{|n|^{4}}$ for $R\geq 1.$\\

\item[(iii)]  $|\Delta_{\md} G(t,n)|\leq \frac{C}{t^{3/2}}$ for $R\leq 1,$ and $|\Delta_{\md} G(t,n)|\leq \frac{C}{|n|^{3}}$ for $R\geq 1.$

\end{itemize}

\end{lemma}
\begin{proof}
(i) First note that by \eqref{prop5} and  $I_{-n}=I_{n}$ for $n\in\N,$ we have that $G(t,n)-G(t,n+1)\geq 0$ and $G(t,-n)-G(t,-n-1)\geq 0$ for $n\in\N_0;$ that is, for each $t\geq 0$ the sequence $G(t,n)$ is increasing for $n\geq 0$ and decreasing for $n\geq 0$ (symmetric). Therefore the equation \eqref{prop2} implies $$|G(t,n)|\leq G(t,0)\leq \frac{C}{t^{1/2}},\quad t>0.$$ On the other hand, we follow some ideas in the proof of \cite[Proposition 3]{C}. Doing a change of variable in the integral given in \eqref{prop4}, we can write $G(t,n)$ in the following way \begin{equation*}\label{G1}G(t,n)=\frac{ 1}{\sqrt{4\pi t}\Gamma(n+1/2)}\int_0^{4t}e^{-u}u^{n-1/2}\biggl(1-\frac{u}{4t}\biggr)^{n-1/2}\,du,\quad n\in\N_0.\end{equation*} Taking into account the inequality \eqref{eq2.2} and identity \eqref{double3}  we get $$G(t,n)\leq C \frac{t}{(n+1)^{3/2}\Gamma(n+1/2)}\int_0^{4t}e^{-u}u^{n-2}\,du\leq \frac{C t\Gamma(n-1)}{n^{3/2}\Gamma(n+1/2)}\leq \frac{C t}{n^3},\quad t>0,n\geq 2.$$ The case $n\leq -1$ follows by symmetry.

(ii) We apply  \eqref{prop5} to get \begin{equation}\label{G2}
G(t,n)-G(t,n+1)=\frac{1}{4t^{3/2}\sqrt{\pi}\Gamma(n+1/2)}\int_0^{4t}e^{-s}s^{n+1/2}\biggl(1-\frac{s}{4t}\biggr)^{n-1/2}\,ds,\quad n\in\N_0.
\end{equation}
On one hand, by previous identity one gets  \begin{eqnarray*}
G(t,n)-G(t,n+1)&\leq&\frac{1}{4t^{3/2}\sqrt{\pi}\Gamma(n+1/2)}\int_0^{4t}e^{-s}s^{n+1/2}\,ds \\
&\leq&\frac{C\Gamma(n+3/2)}{t^{3/2}\Gamma(n+1/2)}\leq \frac{C n} {t^{3/2}},\quad t>0,n\in\N.
\end{eqnarray*}
On the other hand, by equations \eqref{G2}, \eqref{eq2.2} and \eqref{double3} we have
\begin{eqnarray*}
G(t,n)-G(t,n+1)&=&\frac{8t}{\sqrt{\pi}\Gamma(n+1/2)}\int_0^{4t}e^{-s}s^{n-2}\biggl( \frac{s}{4t} \biggr)^{5/2}\biggl(1-\frac{s}{4t}\biggr)^{n-1/2}\,ds \\
&\leq&\frac{C t}{(n+2)^{5/2}\Gamma(n+1/2)}\int_0^{4t}e^{-s}s^{n-2}\,ds \\
&\leq&\frac{C t\Gamma(n-1)}{(n+2)^{5/2}\Gamma(n+1/2)}\leq \frac{C t}{n^4},\quad t>0,n\geq 2.
\end{eqnarray*}

(iii) By \eqref{prop6} we can write  that \begin{eqnarray}\label{G3}
\nonumber &\quad& G(t,n+2)-2G(t,n+1)+G(t,n)\\\nonumber &\quad&\qquad \qquad =\frac{1}{\sqrt{\pi}\Gamma(n+1/2)t^{1/2}}\int_{0}^{4t} e^{-s}\biggl(\frac{s}{2t}-1\biggr)s^{n-1/2}\biggl(1-\frac{s}{4t}\biggr)^{n-1/2}\,\frac{ds}{2t}\\ \nonumber \\
&\quad&\qquad \qquad \qquad \qquad +\frac{1}{\sqrt{\pi}\Gamma(n+1/2)t^{1/2}}\int_{0}^{4t} e^{-s}s^{n+1/2}\biggl(\frac{s}{4t}\biggr)\biggl(1-\frac{s}{4t}\biggr)^{n-1/2}\,\frac{ds}{2t}\\ \nonumber \\
\nonumber &\quad&\qquad \qquad = (I)+(II),
\end{eqnarray}
for $n\in\N_0$. First note that for $n\in\N,$ $|(1-\frac{s}{4t})^{n-1/2}|, |(\frac{s}{2t}-1)| \leq 1$ for $0\leq s\leq 4t.$ Then $$|(I)|\leq \frac{C}{\Gamma(n+1/2)t^{3/2}}\int_0^{4t} e^{-s}s^{n-1/2}\,ds\leq \frac{C}{t^{3/2}},\quad t>0.$$ Secondly, by \eqref{eq2.2} one gets for $n\in\N$ $$|(II)|\leq \frac{C}{(n+1/2)\Gamma(n+1/2)t^{3/2}}\int_0^{4t} e^{-s}s^{n+1/2}\,ds\leq \frac{C}{t^{3/2}},\quad t>0.$$ So, $$|G(t,n+2)-2G(t,n+1)+G(t,n)| \leq  \frac{C}{t^{3/2}},\quad t>0,n\in\N.$$ Furthermore, by \eqref{prop2} we can write $G(t,2)-2G(t,1)+G(t,0)=-\frac{1}{4\sqrt{\pi}t^{3/2}}+O(\frac{1}{t^{5/2}}),$ for $t$ large. Then by symmetry we conclude $|\Delta_\md G(t,n)|\leq \frac{C}{t^{3/2}},$ for $t>0$ and $n\neq 0.$

To show the second inequality,  we consider again \eqref{G3}. Take $n\geq 2.$  Since $| (\frac{s}{2t}-1)| \leq 1$ for $0\leq s\leq 4t,$ we apply \eqref{eq2.2} and \eqref{double3} to get  that \begin{eqnarray*}
|(I)|&\leq& \frac{C}{\Gamma(n+1/2)}\int_0^{4t} e^{-s}s^{n-2}\biggl(\frac{s}{4t} \biggr)^{3/2}\biggl(1-\frac{s}{4t}\biggr)^{n-1/2}\,ds\\
&\leq& \frac{C}{(n+1)^{3/2}\Gamma(n+1/2)}\int_0^{4t} e^{-s}s^{n-2}\,ds\leq \frac{C\Gamma(n-1)}{(n+1)^{3/2}\Gamma(n+1/2)}\leq\frac{C}{(n+1)^3}.
\end{eqnarray*}
Secondly for $n\in\N$, by the same arguments \begin{eqnarray*}
|(II)|&\leq& \frac{C}{\Gamma(n+1/2)}\int_0^{4t} e^{-s}s^{n-1}\biggl(\frac{s}{4t} \biggr)^{5/2}\biggl(1-\frac{s}{4t}\biggr)^{n-1/2}\,ds\\
&\leq& \frac{C}{(n+2)^{5/2}\Gamma(n+1/2)}\int_0^{4t} e^{-s}s^{n-1}\,ds\leq \frac{C\Gamma(n)}{(n+2)^{5/2}\Gamma(n+1/2)}\leq\frac{C}{(n+1)^3}.
\end{eqnarray*} So, by symmetry, we conclude that $$|\Delta_\md G(t,n)|\leq \frac{C}{n^3},\quad t>0,|n|\geq 3.$$
\end{proof}

\begin{remark}\label{Rem2.2}{\it
Observe that the above result is also valid for $n=0$ in the following way. By the proof of the item (i), $G(t,0)\leq \frac{C}{t^{1/2}},$ for $t>0.$ Furthermore by \eqref{prop2} we have $$G(t,0)-G(t,1)=\frac{1}{8\sqrt{\pi} t^{3/2}}+O(\frac{1}{t^{5/2}}),\quad t\text { large,}$$ so $$G(t,0)-G(t,1)\leq \frac{C}{t^{3/2}},\quad t>0.$$ Also, note that $|G(t,1)-2G(t,0)+G(t,-1)|=2|G(t,1)-G(t,0)|\leq \frac{C}{t^{3/2}},$ by the previous comment.

Furthermore, note that for $n\leq -1$ we can write $|\nabla_\md G(t,n)|=|\nabla_\md G(t,-n-1)|,$ and we also have the bounds given in Lemma \ref{Pointwiseestimates}.
}
\end{remark}

 Note we consider the usual gaussian semigroup $(g_t)_{t>0}$ in the Lebesgue spaces $L^p(\R)$. Then it is well-known that
\begin{itemize}
\item[(i)]$\displaystyle{\|g_t\|_{p}\leq \frac{C_p}{t^{{1\over 2}(1-\frac{1}{p})}}}$,\\

\item[(ii)]$\displaystyle{\|g'_t\|_{p}\leq \frac{C_p} {t^{{1\over 2}(1-\frac{1}{p})+\frac{1}{2}}}}$,\\

\item[(iii)] $\displaystyle{ \|g''_t\|_{p}\leq \frac{C_p} {t^{{1\over 2}(1-\frac{1}{p})+1}}},$\\

\end{itemize}
for $t>0.$ To show them, we may check directly $\Vert \quad \Vert_p$ for any $1\le p\le \infty$.  An alternative (and sharp) proof of these inequalities is to check $\Vert \quad \Vert_1$ and $\Vert \quad \Vert_{\infty}$ norms, and apply  Lyapunov inequality:  given $X$  a measure space and $f\in L^1(X) \cap L^\infty(X)$, then
$$
\Vert f\Vert_p\le \Vert f\Vert_1^{1\over p}\Vert f\Vert_\infty^{1-{1\over p}},
$$
for $1<p<\infty$. In the next theorem, we show that the same inequalities hold for the discrete heat semigroup.

\begin{theorem}\label{Kernelestimates} Let $1\leq p\leq \infty.$ Then for $t>0$ we have that
\begin{itemize}
\item[(i)]$\displaystyle{\|G(t,\cdot)\|_p\leq \frac{C_p}{t^{{1\over 2}(1-\frac{1}{p})}}}$.\\

\item[(ii)]$\displaystyle{\|\nabla_\md G(t,\cdot)\|_p\leq \frac{C_p} {t^{{1\over 2}(1-\frac{1}{p})+\frac{1}{2}}}}$.\\

\item[(iii)] $\displaystyle{ \|\Delta_\md G(t,\cdot)\|_p\leq \frac{C_p} {t^{{1\over 2}(1-\frac{1}{p})+1}}}.$
\end{itemize}
\end{theorem}
\begin{proof}
(i) From \eqref{generating} and $I_n(t)\geq 0$ for  $n\in\Z, t\geq0,$ it follows that $$\Vert G(t,\cdot)\Vert_1= \sum_{n\in\Z}|G(t,n)|=1.$$ The proof of item (i) of Lemma \ref{Pointwiseestimates} gives $$\Vert G(t,\cdot)\Vert_\infty= \sup_{n\in\Z}|G(t,n)|=G(t,0)\leq \frac{C}{t^{1/2}}.$$  We apply Lyapunov inequality to get $$\|G(t,\cdot)\|_p\leq \frac{C_p}{t^{{1\over 2}(1-\frac{1}{p})}}.$$

To show (ii), we know that $G(t,n)$ is symmetric on the variable $n\in\Z,$ and the equation \eqref{prop5} implies that it is decreasing on $|n|.$ Since $G(t,\cdot)\in \ell^1,$ then $G(t,|n|)\to0$ as $|n|\to\infty,$ and $$\|\nabla_\md G(t,\cdot)\|_1=\sum_{n\in\Z}|G(t,n+1)-G(t,n)|=2G(t,0)\leq \frac{C}{t^{1/2}}.$$ On the other hand, by equations \eqref{G2}, \eqref{eq2.2} and \eqref{double3} one gets  that $$G(t,n)-G(t,n+1)\leq \frac{C}{tn^{1/2}\Gamma(n+1/2)}\int_0^{4t}e^{-s}s^{n}\leq\frac{C}{t},$$ for $n\in\N$; in Remark \ref{Rem2.2} we have shown that $G(t,0)-G(t,1)\leq \frac{C}{t^{3/2}}\leq \frac{C}{t}$ for $t\geq 1$. The case $t\in (0,1)$ is clear by the continuity of $G(t,0)-G(t,1)$ on $[0,1].$ Therefore $$\|\nabla_\md G(t,\cdot)\|_\infty=\sup_{n\in\Z}|\nabla_\md G(t,n)|\leq \frac{C}{t}.$$ So, by the Lyapunov inequality, we conclude $$\|\nabla_\md G(t,\cdot)\|_p\leq \frac{C_p}{t^{{1\over 2}(1-\frac{1}{p})+\frac{1}{2}}}.$$

Finally, to show (iii), note that $$ \|\Delta_\md G(t,\cdot)\|_1=\sum_{n\in\Z}|\Delta_{\md}G(t,n)|=2\sum_{n=1}^{\infty}|\Delta_{d}G(t,n)|+2(G(t,0)-G(t,1)).$$ On one hand, by previous estimates, $$G(t,0)-G(t,1)\leq \frac{C}{t},\quad t>0.$$ On the other hand, it follows from \eqref{prop6} that \begin{eqnarray*}
G(t,n+2,t)-2G(t,n+1)-G(t,n)&=&\frac{4^{n+1/2}t^{n-1}}{\sqrt{\pi}\Gamma(n+1/2)}\int_{0}^{1} e^{-4ts}s^{n+1/2}(1-s)^{n-1/2}\,ds\\ \\
&&-\frac{4^{n}t^{n-1}}{\sqrt{\pi}\Gamma(n+1/2)}\int_{0}^{1} e^{-4ts}s^{n-1/2}(1-s)^{n-1/2}\,ds\\ \\
&&+\frac{4^{n+1}t^n}{\sqrt{\pi}\Gamma(n+1/2)}\int_{0}^{1} e^{-4ts}s^{n+3/2}(1-s)^{n-1/2}\,ds\\ \\
&=& (I)-(II)+(III),
\end{eqnarray*}
for $n\in\N_0$. Then we apply (\ref{Hyper}) to get
\begin{eqnarray*}
(I)&=&\frac{4^{n+1/2}t^{n-1}e^{-4t}(n+1/2)}{\sqrt{\pi}}\sum_{k=0}^{\infty}\frac{(4t)^k}{k!}\frac{\Gamma(n+1/2+k)}{\Gamma(2n+2+k)}\\
&\leq&\frac{4^{n+1/2}t^{n-1}e^{-4t}}{\sqrt{\pi}}\sum_{k=0}^{\infty}\frac{(4t)^k}{k!}\frac{\Gamma(n+1/2+k)}{\Gamma(2n+1+k)}\\
&=&\frac{4^{n+1/2}e^{-2t}}{t\sqrt{\pi}}\frac{\Gamma(n+1/2)\Gamma(n+1)}{\Gamma(2n+1)}I_n(2t)=C\frac{e^{-2t}}{t}I_n(2t),
\end{eqnarray*}
where  we have used (\ref{lebedev}) in the last line. Also, by \eqref{Hyper} and (\ref{lebedev}) one gets
$$
(II)=\frac{4^{n}e^{-2t}}{t\sqrt{\pi}}\frac{\Gamma(n+1/2)\Gamma(n+1)}{\Gamma(2n+1)}I_n(2t)=C\frac{e^{-2t}}{t}I_n(2t).
$$
Finally, by \eqref{Hyper}, we obtain that \begin{eqnarray*}
(III)&=&\frac{4^{n+1}t^ne^{-4t}(n+3/2)(n+1/2)}{\sqrt{\pi}}\sum_{k=0}^{\infty}\frac{(4t)^k}{k!}\frac{\Gamma(n+1/2+k)}{\Gamma(2n+3+k)}\\
&\leq &\frac{4^{n+4}t^ne^{-4t}(n+3/2)(n+1/2)}{\sqrt{\pi}}\sum_{k=0}^{\infty}\frac{(4t)^k}{k!}\frac{\Gamma(n+5/2+k)}{\Gamma(2n+5+k)},
\end{eqnarray*}
where we have used that $\displaystyle{\frac{2n+3+k}{n+1/2+k},\frac{2n+4+k}{n+3/2+k}\leq 8.}$  We apply again (\ref{lebedev}) to get
\begin{eqnarray*}(III)&\leq& \frac{4^{n+4}e^{-2t}(n+3/2)(n+1/2)}{t^2\sqrt{\pi}}\frac{\Gamma(n+5/2)\Gamma(n+3)}{\Gamma(2n+5)}I_{n+2}(2t)\\
&\leq& C\frac{(n+2)(n+1)e^{-2t}}{t^2}I_{n+2}(2t).\\
\end{eqnarray*}
Since $\sum_{n\in\Z}e^{-2t}I_{n}(2t)=1,$ and, by Corollary \ref{mean}, $\sum_{n\in\Z}n(n-1)e^{-2t}I_{n}(2t)=2t,$ then  the norm $ \|\Delta_\md G(t,\cdot)\|_1$  has the desired bound, that is, $C/t.$

The proof of item (iii) of Lemma \ref{Pointwiseestimates} and Remark \ref{Rem2.2} gives
$$ \|\Delta_\md G(t,\cdot)\|_\infty=\sup_{n\in\Z}|\Delta_\md G(t,n)|\leq\frac{C}{t^{3/2}},$$
and the Lyapunov inequality  concludes the result.\end{proof}

\section{$\ell^p$-$\ell^q$ theory}

Here, we state $\ell^p$-$\ell^q$ decay rates for the solution of \eqref{eq0} assuming certain properties of regularity on the data $f$ and on the forcing term $g$ (Theorem \ref{bounds}). Such assumptions allow $u,$ given in \eqref{sol}, to be a mild solution (see Introduction). Also, we present an optimal $\ell^2$-$\ell^1$ result (Theorem \ref{opti}).

For convenience, we write the solution as $u(t,n)=u_f(t, n) + u_g(t,n)$ where
\begin{eqnarray*}
u_f(t,n)&:=&W_tf(n)=(G(t,\cdot)*f)(n),\\
u_g(t,n)&:=&\int_0^t (W_{t-s}g(s,\cdot))(n)\,ds=\int_0^t (G(t-s,\cdot)*g(s,\cdot))(n)\,ds.
\end{eqnarray*}

\begin{theorem} \label{bounds}
Let $1\le q\le p\le \infty$.

\begin{itemize}
\item[(1)] If $f\in \ell^q$ and $t>0,$ then
\begin{itemize}
 \item[(i)] $\displaystyle{\|u_f(t,\cdot)\|_{p}\leq C t^{-\frac{1}{2}\left(\frac{1}{q}-\frac{1}{p}\right)}\|f\|_{q}}$.\\

 \item[(ii)] $\displaystyle{ \|\nabla_{\md} u_f(t,\cdot)\|_p\leq C t^{-\frac{1}{2}\left(\frac{1}{q}-\frac{1}{p}\right)-\frac{1}{2}}\|f\|_{q},}$\\

 \item[(iii)] $\displaystyle{ \|\Delta_{\md} u_f(t,\cdot)\|_p\leq C t^{-\frac{1}{2}\left(\frac{1}{q}-\frac{1}{p}\right)-1}\|f\|_{q} .}$
\end{itemize}

\medskip

\item[(2)] Let $g(t,\cdot) \in \ell^q,$ and assume there exists $K,\gamma>0$ such that
	  $$\|g(t,\cdot) \|_{q}\leq\frac{K}{(1+t)^{\gamma}},\quad t>0.$$ Then for $t>0$ we have
 \begin{itemize}
\item[(i)] $\|u_g(t,\cdot)\|_{p}\leq  C t^{1-\min\{1,\gamma\}-\frac{1}{2}(\frac{1}{q}-\frac{1}{p})},\quad \gamma\neq 1.$\\

\item[(ii)] $\|u_g(t,\cdot)\|_{p}\leq Ct^{-\frac{1}{2}(\frac{1}{q}-\frac{1}{p})}\log(1+t),\quad \gamma=1.$
\end{itemize}
\end{itemize}
	
\end{theorem}

\begin{proof} (1) We choose $r\ge 1$ such that $1+\dfrac{1}{p}=\dfrac{1}{r}+\dfrac{1}{q}.$ By Young inequality for discrete convolution and Theorem \ref{Kernelestimates} one gets
\begin{align*}
\|u_f(t,\cdot)\|_{p}=\|G(t,\cdot)* f\|_{p}&\le \|G(t,\cdot)\|_{r}\|f\|_{q}\leq C t^{-\frac{1}{2}\left(1-\frac{1}{r}\right)}\|f\|_{q}.
\end{align*}
The estimates for the discrete gradient and Laplacian follow similarly.

(2) Let be again $r\geq1$ satisfying $1+\frac{1}{p}=\frac{1}{q}+\frac{1}{r}$. By Young inequality and Theorem \ref{Kernelestimates} we have \begin{eqnarray*}
\|u_g(t,n)\|_{p}& \leq C &\int_{0}^{t}\frac{1}{(t-s)^{1/2(1-1/r)}}\frac{1}{(1+s)^{\gamma}}\,ds\\
&=& C\biggl(\int_{0}^{t/2}+\int_{t/2}^t\biggr)\frac{1}{(t-s)^{1/2(1-1/r)}}\frac{1}{(1+s)^{\gamma}}\,ds:=I_1+I_2.
\end{eqnarray*}
For the first integral, we have $$I_1\leq \frac{C}{t^{\frac{1}{2}(1-1/r)}}\int_0^{t/2}\frac{1}{(1+s)^{\gamma}}\,ds.$$ On the other hand $$I_2\leq\frac{C}{(1+t)^{\gamma}}\int_0^{t/2}\frac{1}{s^{\frac{1}{2}(1-1/r)}}\,ds\leq C t^{1-\gamma-\frac{1}{2}(1-1/r)}.$$
Finally observe that $I_2$ decays faster than $I_1$ if $\gamma\geq 1,$ and $I_1,I_2$ decay with the same rate if $\gamma\in (0,1).$ So we conclude the result.
\end{proof}

\begin{remark}{\it In order to get bounds for $\displaystyle{ \|\nabla_{\md} u_g(t,\cdot)\|_p}$ and $\displaystyle{ \|\Delta_{\md} u_g(t,\cdot)\|_p}$ in Theorem \ref{bounds} (2), we would have to impose certain extra regularity conditions to $g.$}
\end{remark}

Next result shows that the previous  $\ell^2$-$\ell^1$ decay is optimal. We include an alternative proof of the upper bound to get a self-contained theorem. The proof is based on standard Fourier techniques in the Hilbert $\ell^2$ space .

\begin{theorem}\label{opti}
Assume that $f\in \ell^1$ and $\sum_{n\in\Z}f(n)\neq 0.$ Then for large enough $t$ there exist $c,C>0$ such that
$$
\frac{c}{t^{1/4}}\biggl| \sum_{n\in\Z}f(n) \biggr|\leq \|u_f(t,\cdot)\|_2 \leq \frac{C}{t^{1/4}}\|f\|_1.
$$
\end{theorem}

\begin{proof}
Let $\varepsilon>0$ small enough, then
\begin{align}\label{5.1}
\nonumber\|u_f(t,\cdot)\|_2^2&=\|\mathcal{F}u_f(t,\cdot)\|_2^2=\frac{1}{2\pi}\int_{\mathbb{T}}|\mathcal{F}G(t,\cdot)(\theta)|^{2}\,|\mathcal{F}f(\theta)|^2\,d\theta\\
&\geq \frac{1}{2\pi}\int_{-\varepsilon}^{\varepsilon}e^{-8t\sin^2(\theta/2)}\,\,|\mathcal{F}f(\theta)|^2\,d\theta\\
\nonumber&\geq \frac{e^{-8t\sin^2(\varepsilon/2)}}{2\pi}\int_{-\varepsilon}^{\varepsilon}|\mathcal{F}f(\theta)|^2\,d\theta.
\end{align}
Since $f\in\ell^1$ then $\mathcal{F}f\in C(\mathbb{T})\cap L^2(\T)$. By the Lebesgue differentiation theorem, we may choose $\varepsilon_0$ small enough such that
$$
\frac{1}{2\varepsilon}\int_{-\varepsilon}^{\varepsilon}|\mathcal{F}f(\theta)|^2\,d\theta\geq \frac{1}{2}\left|\mathcal{F}f(0)\right|^2\quad
\mbox{for all}\quad \varepsilon\in(0,\varepsilon_0].
$$

Substituting the previous inequality in \eqref{5.1} we have that for all $\varepsilon\in(0,\varepsilon_0]$
$$
\|u_f(t,\cdot)\|_2^2\geq \frac{\varepsilon e^{-8t\sin^2(\varepsilon/2)}}{2\pi}\left|\mathcal{F}f(0)\right|^2.
$$
We choose $\varepsilon:=\dfrac{\varepsilon_0}{t^{1/2}}<\varepsilon_0$ for large enough $t$. Note that for large $t$ we have $8t\sin^2(\frac{\varepsilon_0}{2t^{1/2}})$ is bounded. Therefore $$
\|u_f(t,\cdot)\|_2^2\geq \frac{\varepsilon e^{-8t\sin^2(\varepsilon/2)}}{2}\left|\mathcal{F}f(0)\right|^2=\frac{\varepsilon_0 e^{-8t\sin^2(\frac{\varepsilon_0}{2t^{1/2}})}}{2t^{1/2}}\left|\mathcal{F}f(0)\right|^2\geq \frac{c^2}{t^{1/2}}\biggl| \sum_{n\in\Z}f(n) \biggr|^2,$$
with $c$ a positive constant.

Next, let us prove the upper bound. By Plancherel's Theorem, we have

\begin{align*}
\nonumber\|u_f(t,\cdot)\|_2^2&=\|\mathcal{F}u_f(t,\cdot)\|_2^2=\frac{1}{2\pi}\int_{\mathbb{T}}|\mathcal{F}G(t,\cdot)(\theta)|^{2}\,|\mathcal{F}f(\theta)|^2\,d\theta\\
\nonumber&= \frac{1}{2\pi}\int_{\mathbb{T}}e^{-8t\sin^2(\frac{\theta}{2})}\,\,|\mathcal{F}f(\theta)|^2\,d\theta\leq\frac{1}{2\pi}\|\mathcal{F}f\|_{\infty}^2 \int_{\mathbb{T}}e^{-8t\sin^2(\frac{\theta}{2})}\,d\theta\\
&\leq\frac{1}{2\pi}\|f\|_{1}^2 \int_{\mathbb{T}}e^{-8t\sin^2(\frac{\theta}{2})}\,d\theta.\\
\end{align*}
Note that by the inverse Fourier transform we can write $$G(2t,0)=\frac{1}{2\pi}\int_{\T}e^{-8t\sin^2(\frac{\theta}{2})}\,d\theta.$$ Then by \eqref{prop2} we have $$\nonumber\|u_f(t,\cdot)\|_2^2\leq \frac{C^2}{t^{1/2}}\|f\|_1^2,$$ with $C>0.$
\end{proof}


\section{Large time behaviour of solutions}\label{Non-homogeneous Problem}

In this section we study the large time behaviour of the solution of \eqref{eq0}. We prove that the solution $u=u_f+u_g$ converges asymptotically on $\ell^p$ to a constant (given by the sum of the mass of the initial data $f$ and the mass of the linear forcing term $g$) times the discrete heat kernel (Theorem \ref{ultimo}). In particular, we get the rate of the convergence. Along the section we will assume the following:

\begin{enumerate}
\item[(a)] $f\in \ell^1$.
\item[(b)] There exist $K>0,$ $\gamma>1$  such that
$$
\|g(t,\cdot)\|_1 \leq \frac{K}{(1+t)^{\gamma}},\quad t>0.
$$
\end{enumerate}
Set $$M_f :=\sum_{n\in\Z} f(n), \quad M_g:=\int_{0}^{\infty}\sum_{n\in\Z} g(t,n)\,dt.$$ We will also assume that $M_f\neq 0,$ since in other case ($M_f=0$), by the mass conservation principle, the solution $u_f$ would be null.

\begin{theorem}\label{ultimo}
	Let $1\leq p\leq \infty.$ Assume that conditions $(a)$ and $(b)$ hold, and let $u$ be the solution of \eqref{eq0}.
	\begin{itemize}
		
		\item[(i)] Then $$t^{\frac{1}{2}(1-\frac{1}{p})}\|u_{f}(t,\cdot)-M_f G(t,\cdot)\|_{p}\to 0,\quad \mbox{as}\quad t\to\infty,$$ and $$t^{\frac{1}{2}(1-\frac{1}{p})}\|u_{g}(t,\cdot)-M_{g}G(t,\cdot)\|_{p}\to 0,\quad \mbox{as}\quad t\to\infty.$$
		
		\item[(ii)] Suppose in addition that $nf\in \ell^1,$ then $$ t^{\frac{1}{2}(1-\frac{1}{p})}\|u_f(t,\cdot)-M_{f}G(t,\cdot)\|_{p}=O\biggl( \frac{1}{t^{1/2}} \biggr),\quad t \text{ large}.$$
		
	\end{itemize}
	
\end{theorem}

\begin{proof}
For convenience in some particular points of the proof and without loss of generality, we prove the result for non-negative sequences, since in general case we would write $f=f_+-f_-,$ where \begin{equation*}
f_+(n)=\left\{\begin{array}{ll}
f(n),&f(n)>0,\\ \\
0,& f(n)\leq 0,
\end{array} \right.\text{ and }\quad
f_-(n)=\left\{\begin{array}{ll}
-f(n),&f(n)<0,\\ \\
0,& f(n)\geq 0,
\end{array} \right.
\end{equation*}
both are non-negative sequences.

We start with assertion $(ii)$.  We write \begin{eqnarray*}
u_f(t,n)-M_fG(t,n)&=&\sum_{j\leq -1} (G(t,n-j)-G(t,n))f(j)+\sum_{j\geq 1} (G(t,n-j)-G(t,n))f(j)\\
&=&\sum_{j\leq -1}\sum_{k=j+1}^{0}\nabla_\md G(t,n-k)f(j)-\sum_{j\geq 1} \sum_{k=0}^{j-1}\nabla_\md G(t,n-k-1)f(j)\\
&=&\sum_{k=-\infty}^0 \nabla_\md G(t,n-k)\sum_{j=-\infty}^{k-1}f(j)-\sum_{k=1}^{\infty} \nabla_\md G(t,n-k)\sum_{j=k}^{\infty}f(j)\\
&=& (\nabla_\md G(t,\cdot)*F)(n),
\end{eqnarray*}
where $F(k)=\sum_{j=-\infty}^{k-1}f(j)$ for $k\leq 0,$ and $F(k)=-\sum_{j=k}^{\infty}f(j)$ for $k\geq 1.$ Note that $$
\|F\|_1\leq\sum_{k=-\infty}^{0}\sum_{j=-\infty}^{k-1}|f(j)|+\sum_{k=1}^{\infty}\sum_{j=k}^{\infty}|f(j)|=\sum_{j\neq 0} |jf(j)|.$$
Then by Theorem \ref{Kernelestimates} it follows that  $$\|u_f(t,\cdot)-M_fG(\cdot,t)\|_{p}\leq \frac{C\|F\|_1}{t^{\frac{1}{2}(1-\frac{1}{p})+\frac{1}{2}}},\quad t>0.$$

Now we focus on the first part of assertion (i). We take the sequence $(f_{j})_{j\in\N_0}$ given by

$$
f_{j}(n)=\left\{\begin{array}{ll}
\displaystyle{\frac{M_f}{\displaystyle\sum_{|k|\leq j}f(k)}f(n)},&|n|\leq j,\\ \\
0,& |n|>j,
\end{array} \right.
$$
whenever $f(0)\neq 0.$ Else, if $j_0$ is the first natural number such that $f(j_0)$ and/or $f(-j_0)$ are positive, then we consider the sequence $(f_{j})_{j\geq j_0}$. It is an easy computation that $$\displaystyle\sum_{n\in\Z} f_{j}(n)=M_f,\text{ for all }j\in\N_0,$$ and $$\|f_j-f\|_{1}\to 0,\quad j\to \infty.$$ We write $$\|u_f(t,\cdot)-M_fG(t,\cdot)\|_{p}\le \|G(t,\cdot)\|_p \|f_j-f\|_{1}+\|G(t,\cdot)* f_j-M_fG(t,\cdot)\|_{p}.$$ Note that $(nf_j(n))$ belongs to $\ell^1$ for each $j.$ Then, by the proof of the statement $(ii)$ and Theorem \ref{Kernelestimates}, it follows that
$$t^{\frac{1}{2}\left(1-\frac{1}{p}\right)}\|u_f(t,\cdot)-M_fG(t,\cdot)\|_{p}\leq  C\|f_j-f\|_{1} +C_jt^{-\frac{1}{2}},\quad t>0,$$
where $C_j$ depends on $\|nf_j(\cdot)\|_{1}.$ Therefore, $$\limsup_{t\to\infty}t^{\frac{1}{2}\left(1-\frac{1}{p}\right)}\|u_f(t,\cdot)-M_fG(t,\cdot)\|_{p}\leq  C\|f_j-f\|_{1} ,$$ and the result follows by taking $j\to \infty.$

Next, let us prove the second part of $(i)$. We can write
	
	\begin{equation*}
	M_{g}=\int_{0}^{t}\sum_{m\in\Z}^{} g(s,m)\,ds+\int_{t}^{\infty}\sum_{m\in\Z}^{} g(s,m)\,ds,\quad t>0.
	\end{equation*}
	It follows from Theorem \ref{Kernelestimates} that
	\begin{align*}
	&t^{\frac{1}{2}(1-\frac{1}{p})}\left\|G(t,\cdot)\int_{t}^{\infty}\sum_{m\in\Z}^{}g(s,m)\,ds\right\|_{p}\\
	&\leq t^{\frac{1}{2}(1-\frac{1}{p})}\|G(t,\cdot)\|_{p}
	\int_{t}^{\infty}\sum_{m\in\Z}^{} |g(s,m)|\,ds\to 0,\quad t\to\infty.
	\end{align*}
	Therefore in the following we will prove that
	\begin{equation*}
	t^{\frac{1}{2}(1-\frac{1}{p})}\left\|\int_{0}^{t}(G(t-s,\cdot)\ast g(s,\cdot))(\cdot)-G(t,\cdot)\int_{0}^{t}\sum_{m\in\Z}^{} g(s,m)\,ds\right\|_{p}\to 0,\quad t\to\infty.
	\end{equation*}
	
We take $0<\delta<\frac{1}{10},$ which implies in particular $0<\delta<\dfrac{1}{5}<\dfrac{1}{2}$ and
	\begin{equation}\label{delta}
	\frac{\delta}{1-\delta}<\frac{1}{4}.
	\end{equation}
Let $t>0$ fixed. We decompose the set $(0,t)\times\Z$ into two parts
	$$\Omega_1= (0,\delta t)\times\{m\in\mathbb{Z}:|m|\leq \lfloor(\delta t)^{\frac{1}{2}}\rfloor\},\quad \Omega_2=((0,t)\times\Z)\setminus \Omega_1,$$ where $\lfloor (\delta t)^{\frac{1}{2}}\rfloor$ denotes the integer part of the positive number $(\delta t)^{\frac{1}{2}}.$

Let us start with the set $\Omega_1,$ and we write
	\begin{align*}
	&\left\|\int\sum_{(s,m)\in \Omega_1}(G(t-s,\cdot-m)-G(t,\cdot))g(s,m)\ ds\right\|_{p}\\
	&\leq\int \sum_{(s,m)\in \Omega_1}\| G(t-s,\cdot-m)-G(t,\cdot)\|_{p}|g(s,m)|\ ds.
	\end{align*} Note that in this set the following inequalities hold
	\begin{equation}\label{inequalities}
	t\geq t-s\geq t(1-\delta)>\frac{t}{2}.
	\end{equation}
For each $(s,m)\in\Omega_1,$ we consider the following sets on $\Z,$ $$A= \{n\in\mathbb{Z}:|n-m|<2\lfloor(\delta t)^{\frac{1}{2}}\rfloor\},\quad B= \{n\in\mathbb{Z}:|n-m|\geq 2\lfloor (\delta t)^{\frac{1}{2}}\rfloor\},$$ and we write the $p$-norm over $\Omega_1$ in the following way \begin{align*}
	\| G(t-s,\cdot-m)-G(t,\cdot)\|_{p}&\leq
	\left(\sum_{n\in A}|G(t-s,n-m)-G(t,n)|^p\right)^{\frac{1}{p}}\\
	&+\left(\sum_{n\in B}|G(t-s,n-m)-G(t,n)|^p\right)^{\frac{1}{p}}.\\
	\end{align*}
Now we focus on $A.$ We write $$\biggl(\sum_{n\in A}|G(t-s,n-m)-G(t,n)|^p\biggr)^{1/p}\leq \biggl(\sum_{n\in A}|G(t-s,n-m)|^p\biggr)^{1/p}+\biggl(\sum_{n\in A}|G(t,n)|^p\biggr)^{1/p}:=I_1+I_2.$$

For $(s,m)\in \Omega_1$ and $n\in A$ we have by \eqref{delta} and \eqref{inequalities} that $\frac{|n-m|^2}{t-s} < \frac{(2\lfloor (\delta t)^{1/2}\rfloor)^2}{t(1-\delta)}\leq\frac{4\delta}{(1-\delta)}<1,$ so it follows from Lemma \ref{Pointwiseestimates} that $$ t^{\frac{1}{2}(1-\frac{1}{p})}I_1\leq C \frac{t^{\frac{1}{2}(1-\frac{1}{p})}}{(t-s)^{\frac{1}{2}}}
	\left(\sum_{n\in A}1\right)^{\frac{1}{p}}\leq  C\frac{t^{\frac{1}{2}(1-\frac{1}{p})}}{t^{1/2}(1-\delta)^{1/2}}(4\lfloor(\delta t)^{\frac{1}{2}}\rfloor)^{1/p}\leq  C \delta^{\frac{1}{2p}}.$$ Analogously, we compute
	$$\frac{|n|^2}{t}\leq \frac{(|n-m|+|m|)^2}{t}\leq 9\delta<1,$$ and therefore by Lemma \ref{Pointwiseestimates} we get
	\begin{align*}
	t^{\frac{1}{2}(1-\frac{1}{p})}I_2\leq  C \delta^{\frac{1}{2p}}.
	\end{align*}
	Since $\int_0^\infty \|g(s,\cdot)\|_1\,ds<\infty$ we have
	$$
	t^{\frac{1}{2}(1-\frac{1}{p})}\int\sum_{(s,m)\in \Omega_1}\left(\sum_{n\in A}^{}|G(t-s,n-m)-G(t,n)|^p\right)^{\frac{1}{p}} |g(s,m)|\ ds \leq C \delta^{\frac{1}{2p}}\to 0,\quad \delta\to 0.
	$$
	
Next, we consider on $\Omega_1$ the part $B$ of the $p$-norm. Let us write \begin{align*}
	\left(\sum_{n\in B}|G(t-s,n-m)-G(t,n)|^p\right)^{\frac{1}{p}}&\leq \left(\sum_{n\in B}|G(t-s,n-m)-G(t-s,n)|^p\right)^{\frac{1}{p}}\\
	&+\left(\sum_{n\in B}|G(t-s,n)-G(t,n)|^p\right)^{\frac{1}{p}}\\
	&=:I_3+I_4.
\end{align*}

First, we focus on estimate $I_3$. Observe that if $m\geq 1$ there exists $\tilde{n}\in [n-m,n-1]$ such that $$|G(t-s,n-m)-G(t-s,n)|\leq |m| |\nabla_\md G(t-s,\tilde{n})|,$$ and if $m\leq -1$ there exists $\tilde{n}\in [n,n-m-1]$ such that $$|G(t-s,n-m)-G(t-s,n)|\leq |m| |\nabla_\md G(t-s,\tilde{n})|.$$ Therefore $$I_3\leq |m|\left(\sum_{n\in B}|\nabla_\md G(t-s,\tilde{n})|^p\right)^{\frac{1}{p}}.$$ The case $m=0$ is trivial. Since $|m|\leq \lfloor(\delta t)^{\frac{1}{2}}\rfloor\leq \frac{1}{2}|n-m|$, then we get
	\begin{equation}\label{comparison1}
	|\tilde{n}|\geq|n-m|-|\tilde{n}-(n-m)|\geq |n-m|-(|m|-1)\geq \frac{|n-m|}{2},
	\end{equation}
	and
	\begin{equation}\label{comparison2}
	|\tilde{n}|\leq|n-m|+(|m|-1)\leq |n-m|+\frac{|n-m|}{2}\leq \frac{3|n-m|}{2},\quad m\geq 1,
	\end{equation}
and \begin{equation}\label{comparison3}
	|\tilde{n}|\leq|n-m-1|+(|m|-1)\leq |n-m-1|+\frac{|n-m|}{2}\leq \frac{5|n-m|}{2},\quad m\leq -1,
	\end{equation}
where in the last inequality we have used that $|n-m-1|\leq 2|n-m|$ since we can assume $n\neq m$ in $B$ (note that if $n=m,$ then $\frac{|n-m|}{2}\geq |m|$ implies $m=0,$ and therefore $n=0.$ This case is trivial since $|G(t-s,n-m)-G(t-s,n)|$ would be null). So, equations \eqref{comparison1}, \eqref{comparison2} and \eqref{comparison3} show that $|\tilde{n}|$ and $|n-m|$ are comparable.

We divide $I_3$ into two parts $I_{31}$ and $I_{32},$ depending on whether $\frac{|\tilde{n}|}{(t-s)^{\frac{1}{2}}}$ is less or greater than 1. By \eqref{comparison1}, the set of summation index on $I_{31}$ is contained in $\{n\in \Z\,:\,1\leq |n-m|\leq 2(t-s)^{1/2}\}.$ Then Lemma \ref{Pointwiseestimates} and \eqref{inequalities} imply \begin{align*}I_{31}&\leq C (\delta t)^{\frac{1}{2}}\left(\sum_{1\leq |n-m|\leq 2(t-s)^{1/2}}\frac{|n-m|^p}{(t-s)^{\frac{3p}{2}}}\right)^{1/p}\leq C\frac{(\delta t)^{\frac{1}{2}}}{t^{\frac{3}{2}}} \left(\sum_{1\leq |n-m|\leq 2(t-s)^{1/2}}|n-m|^p\right)^{1/p}\\ \\
&\leq C\delta^{\frac{1}{2}} t^{-1} (t-s)^{(p+1)/2p}\leq C \delta^{\frac{1}{2}} t^{-\frac{1}{2}(1-\frac{1}{p})}.
\end{align*}
Now, by \eqref{comparison2} and \eqref{comparison3}, the set of summation index on $I_{32}$ is contained in $\{n\in \Z\,:\,|n-m|\geq \frac{2}{5}(t-s)^{1/2}\}.$ Then Lemma \ref{Pointwiseestimates} and \eqref{inequalities} imply
\begin{align*}
I_{32}&\leq C(\delta t)^{\frac{1}{2}}\left(\sum_{|n-m|\geq \frac{2}{5}(t-s)^{1/2}}\left(\frac{t-s}{|n-m|^4}\right)^p\right)^{\frac{1}{p}}\leq C\delta^{\frac{1}{2}}t^{3/2}\left(\sum_{|n-m|\geq \frac{2}{5}(t-s)^{1/2}}\frac{1}{|n-m|^{4p}}\right)^{\frac{1}{p}}\\ \\
&\leq C\frac{\delta^{\frac{1}{2}}t^{3/2}}{(t-s)^{(4p-1)/2p}}\leq C \delta^{\frac{1}{2}} t^{-\frac{1}{2}(1-\frac{1}{p})}.
	\end{align*}
So we get that $t^{\frac{1}{2}(1-\frac{1}{p})}I_3\leq C \delta^{\frac{1}{2}}.$

Next, let us estimate $I_4$.  From  the mean value theorem, there exists $\tilde{t}\in(s,t)$ such that \begin{equation*}
	I_4=s\left(\sum_{n\in B}|\partial_tG(\tilde{t},n)|^p\right)^{\frac{1}{p}}=s\left(\sum_{n\in B}|\Delta_\md G(\tilde{t},n)|^p\right)^{\frac{1}{p}}.
\end{equation*}
Note that in $\Omega_1$ we have $0<s<\delta t.$ On the other hand, \eqref{inequalities} implies  $\frac{t}{2}<t-s\leq \tilde{t}\leq t,$ that is, $\tilde{t}$ and $t$ are comparable. Again, we divide $I_4$ into two cases. We denote by $I_{41}$ and $I_{42}$ depending on whether $\frac{n}{\tilde{t}^{1/2}}$ is less or greater than 1.
	
By previous comments, the set of summation index on $I_{41}$ is contained in  $\{n\in \Z: |n|\leq t^{\frac{1}{2}}\}.$ Then by Lemma \ref{Pointwiseestimates} one gets \begin{align*}
	I_{41}&\leq C \delta t \left(\sum_{ |n|\leq t^{\frac{1}{2}}}\frac{1}{t^{\frac{3p}{2}}}\right)^{\frac{1}{p}}\leq C\delta t^{-\frac{1}{2}(1-\frac{1}{p})}.
	\end{align*}

Next we estimate $I_{42}$. Now the set of summation index on $I_{41}$ is contained in  $\{n\in \Z: |n|\geq \frac{t^{1/2}}{\sqrt{2}}\}.$ It follows by Lemma \ref{Pointwiseestimates} that \begin{align*}
    I_{42}\leq C\delta t \left(\sum_{|n|\geq \frac{t^{1/2}}{\sqrt{2}}}\frac{1}{|n|^{3p}}\right)^{\frac{1}{p}}\leq \frac{C \delta t}{t^{(3p-1)/2p}}=C\delta t^{-\frac{1}{2}(1-\frac{1}{p})}.
    \end{align*}
Therefore $t^{\frac{1}{2}(1-\frac{1}{p})}I_4\leq C\delta.$	

Collecting all above terms over $B,$ there is $\eta>0$ such that
	\begin{align*}\label{generalestimate2}
	&t^{\frac{1}{2}(1-\frac{1}{p})}\int\sum_{(s,m)\in \Omega_1} \left(\sum_{n\in B}|G(t-s,n-m)-G(t,n)|^p\right)^{\frac{1}{p}}|g(s,m)|\,ds \\
	\nonumber &\leq C \delta^{\eta}\int_{0}^{\infty}\sum_{n\in\Z}|g(s,m)|\,ds\to 0,\quad \delta\to 0.
	\end{align*}
Note that the above convergence to zero as $\delta\to 0$ is uniformly in $t$.

Now, we focus on $\Omega_2$. Then \begin{align*}
&t^{\frac{1}{2}(1-\frac{1}{p})}\int\sum_{(s,m)\in \Omega_2}\| G(t-s,\cdot-m)-G(t,\cdot)\|_{p}|g(s,m)|\,ds\\
&\leq t^{\frac{1}{2}(1-\frac{1}{p})} \int\sum_{(s,m)\in \Omega_2}\| G(t-s,\cdot-m)\|_{p}|g(s,m)|\,ds\\
&+ t^{\frac{1}{2}(1-\frac{1}{p})} \int\sum_{(s,m)\in \Omega_2}\| G(t,\cdot)\|_{p}|g(s,m)|\,ds\\
&:=I_5+I_6.
	\end{align*}
It follows by Theorem \ref{Kernelestimates} that $$ I_6\leq C \int\sum_{(s,m)\in \Omega_2}|g(s,m)|\,ds.$$ Observe that the set $\Omega_1\to (0,\infty)\times\Z $ as $t\to\infty$, so $\Omega_2$ will have null measure in such case. Since $\int_{0}^{\infty}\sum_{n\in\Z}^{}|g(s,n)|\,ds<\infty$ we conclude $$\int\sum_{(s,m)\in \Omega_2}|g(s,m)|\,ds\to 0, \quad t\to\infty,$$ and therefore $I_6$ also vanishes as $t\to\infty.$
	
Finally, to estimate $I_5,$ we consider the following two cases; either $s< \delta t$ or $s\geq \delta t$. That is, we write
\begin{align*}
I_5&= t^{\frac{1}{2}(1-\frac{1}{p})}\int_{0}^{\delta t}\sum_{|m|>\lfloor(\delta t)^{\frac{1}{2}}\rfloor} \|G(t-s,\cdot-m)\|_{p}|g(s,m)|\,ds\\
&+t^{\frac{1}{2}(1-\frac{1}{p})}\int_{\delta t}^{t}\sum_{m\in \Z}^{} \|G(t-s,\cdot-m)\|_{p}|g(s,m)|\,ds\\
&=:I_{51}+I_{52}.
	\end{align*}
By Theorem \ref{Kernelestimates} and \eqref{inequalities} we have
\begin{align*}
I_{51}&\leq Ct^{\frac{1}{2}(1-\frac{1}{p})}\int_{0}^{\delta t}\sum_{|m|>\lfloor(\delta t)^{\frac{1}{2}}\rfloor} \frac{1}{(t-s)^{\frac{1}{2}(1-\frac{1}{p})}}|g(s,m)|\,ds\\
&\leq C\int_{0}^{\delta t}\sum_{|m|>\lfloor(\delta t)^{\frac{1}{2}}\rfloor} |g(s,m)|\,ds \to 0,\quad t\to\infty.
\end{align*}
Also by Theorem \ref{Kernelestimates} and hypothesis (b) at the beginning of this section, we obtain
$$
I_{52}\leq C t^{\frac{1}{2}(1-\frac{1}{p})}\int_{\delta t}^{t}\sum_{m\in\Z}^{}\dfrac{|g(s,m)|}{(t-s)^{\frac{1}{2}(1-\frac{1}{p})}}\,ds\leq C t^{\frac{1}{2}(1-\frac{1}{p})-\gamma}\int_{\delta t}^{t}\dfrac{1}{(t-s)^{\frac{1}{2}(1-\frac{1}{p})}}\,ds\leq C t^{1-\gamma}\to 0,
	$$
as $t\to\infty,$ and we conclude the proof.
\end{proof}

\begin{remark}{\it Similarly, in the continuous case, it is known that the solution $v$ of \eqref{eq6} satisfies $$t^{\frac{1}{2}(1-\frac{1}{p})}\|v(t,\cdot)-M g_t\|_{p}\to 0,\quad t\to\infty,$$ if $f\in L^1(\R),$ and $$ t^{\frac{1}{2}(1-\frac{1}{p})}\|v(t,\cdot)-M g_t\|_{p}=O\biggl( \frac{1}{t^{1/2}} \biggr),\quad t \text{ large},$$
      if $(1+\vert s\vert)f\in L^1(\R)$, where in both cases $M=\int_{-\infty}^\infty f(s)\,ds$ (see \cite[Th\'eor\`eme 4]{DZ}, \cite[Lemma 2 and 3]{ZE} and in detail \cite[p. 15]{Z}).
    }
\end{remark}



\section{Graphics, open questions, and final comments}

In this last section, we present additional information about three different topics which have appeared along the paper. In fact, the three subsections contain some open problems or future research lines which could be treated in forthcoming works.

\subsection{The sequence of polynomials $(p_k)_{k\ge 0}$}

\begin{figure}[h]
\caption{}
\medskip
\centering
\includegraphics[width=0.8\textwidth]{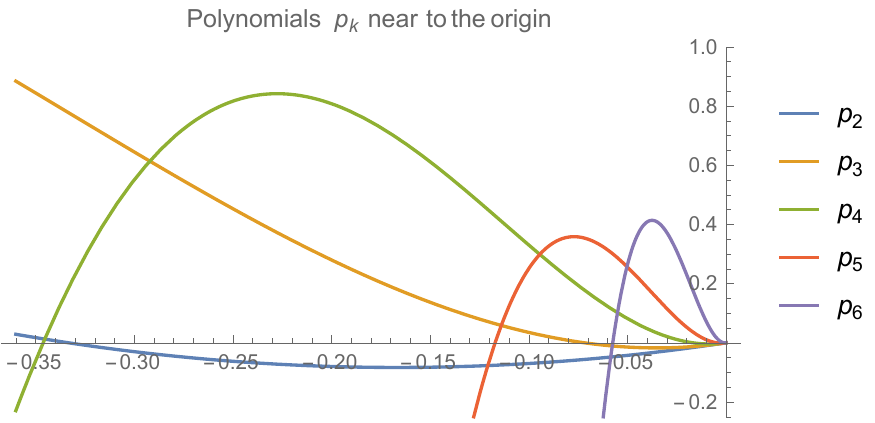}\label{Polinomios}
\end{figure}

We consider the sequence of polynomials $(p_k)_{k\ge 0}$ given in Lemma \ref{sequence}. Note that $p_k(t)> 0$ for $t>0$. In Table \ref{table1}, we present approximations of the zeros of $p_k$ for $2\le k\le 6$. We show these polynomials for close values to the origin in Figure \ref{Polinomios}.

Note that $p_k$ has $k$-non positive real zeros. Moreover, these zeros have the following interlacing property: if $k<j$, there is a zero of $p_j$ between any two zeros of $p_k$. The set of zeros of $p_k$ is the set of zeros of the moments of order $2k$ of Bessel functions $(I_n)_{n\ge 0}$ (Theorem \ref{pepe}). It would be interesting to find a physical interpretation of such zeros. Maybe, the well-standing theory of orthogonal polynomials may be applied to the sequence  $(p_k)_{k\ge 0}$ in order to define a positive measure on $(-\infty, 0)$ (perhaps in terms of Bessel functions), and obtain additional information about these special functions.

\begin{table}[htb]
	\centering
\caption{Zeros of $p_k$ for $ 2\le k\le 6$. }
\begin{tabular}{|l|l|l|l|l|l|l|}
\hline
$p_k$ & $t_0$ &$t_1$ &$t_2$ &$t_3$ &$t_4$&$t_5$\\
\hline
$p_2$ &$0$& $ -0.3333$ & & & &\\
\hline
$p_3$ &$0$& $ -0.0718$ & $ -0.9281$ & & &\\
\hline
$p_4$&$0$& $ -0.01680\dots$ & $ -0.34615\dots$ & $ -1.63703\dots$ & &\\
\hline
$p_5$&$0$& $ -0.00465\dots$ & $ -0.11705\dots$ & $ -0.80786\dots$ &$ -2.40376\dots $&\\
\hline
$p_6$&$0$& $ -0.00099\dots $ & $-0.05743\dots $ & $-0.39502\dots $ &$  -1.31451\dots $&$\: -3.23203\dots$
\\
\hline
\end{tabular}\label{table1}
\end{table}

\subsection{Decayment of finite differences of arbitrary order}
The gaussian semigroup $(g_t)_{t>0}$ in the Lebesgue spaces $L^p(\R)$ satisfies
$$
\displaystyle{ \|g^{(n)}_t\|_p\leq {C_p}{t^{-\frac{1}{2}(1-\frac{1}{p})-{n\over 2 }}},}$$
for $n\ge 0$, see for example \cite[Section 2.2]{Z}. In Theorem \ref{Kernelestimates}, we have proved estimates for finite difference operators of order 1 and 2 applied to the discrete heat semigroup (forward difference and the discrete Laplacian). It seems natural to conjecture that any finite difference operator of order $n$ applied to the discrete heat semigroup has the same decay rate than the $n$-derivative of the continuous heat semigroup.

\subsection{Self similar form for the heat semigroup}
The continuous  heat semigroup has the following self similar form
$$
g_t(s)= {1\over \sqrt{t}}g_1({s\over \sqrt{t}}), \qquad s \in \R,\,\, t>0.
$$
This property allows to show, in a direct way, decay estimates, see for example \cite{DZ,ZE, Z}. However the discrete heat semigroup $(G(t, \cdot))_{t>0}$ seems not to have such property.

Note that, if we change the variable $t=e^{2r}$ for $r\in \R$,  we get
$$
g_t(s)= e^{-r}g_1(e^{-r}s)=:T_{r}(g_1)(s), \qquad s\in \R.
$$
The family $(T_r)_{r\in \R}$ is a one-parameter $C_0$-group of isometries in $L^1(\R)$. It seems natural to ask about a similar $C_0$-group of isometries in $\ell^1$ to attack the problem of self similar form in $\ell^1$ from a different point of view.\\


\end{document}